\documentclass[oneside,a4paper]{article}
\usepackage{akp1}
\usepackage{hyper}
\usepackage{letters}
\usepackage{notation}
\usepackage{hyphenation}
\usepackage{envmnt}

\hypersetup{colorlinks=true,citecolor=black,linkcolor=black,anchorcolor=black,filecolor=black,menucolor=black,urlcolor=black,pdftitle={Alexandrov meets Kirszbraun},pdfauthor={S. Alexander, V. Kapovitch, A. Petrunin},pdfdisplaydoctitle}

\newcommand*\chapter[1]{}
\def\claim#1{%
\par%
\medskip%
\noindent%
\refstepcounter{thm}%
\hbox{\bf\boldmath \arabic{section}.\arabic{thm}. #1.}
\it\ 
}
\def\endclaim{
\par
\medskip}

\begin{document}
\title{Alexandrov meets Kirszbraun}
\author{S. Alexander, V. Kapovitch, A. Petrunin}
\date{}
\maketitle

\begin{abstract}
We give a simplified proof of the generalized Kirszbraun theorem for Alexandrov spaces,
which is due to Lang and Schroeder.
We also discuss related questions, both solved and open. 
\end{abstract}

\section{Introduction}

Kirszbraun's theorem states that any \emph{short map} (i.e. 1-Lipschitz map) from a subset of Euclidean space to another in Euclidean space can be extended as a short map to the whole space.

This theorem was proved first by Kirszbraun in \cite{kirszbraun}.  Later it was reproved by Valentine in \cite{valentine-sphere} and \cite{valentine-kirszbraun}, where he also generalized it to pairs of Hilbert spaces of arbitrary dimension 
as well as pairs of spheres of the same dimension 
and pairs of hyperbolic spaces with the same curvature.
J.~Isbel in \cite{isbell} studied target spaces that satisfy the above condition for any source space.

Valentine was also interested in pairs of metric spaces,
say $\spc{U}$ and $\spc{L}$,
which satisfy the above property, namely, given a subset $Q\subset\spc{U}$, 
any short map $Q\to\spc{L}$ 
can be extended to a short map $\spc{U}\to \spc{L}$.
It turns out that this property has a lot in common with the definition of Alexandrov spaces
(see theorems \ref{thm:kirsz+}, \ref{thm:kirsz-def} and \ref{thm:cba-kirsz-def}).
Surprisingly, this relationship was first discovered only in the 1990's;
it was first published by Lang and Schroeder in \cite{lang-schroeder}.
(The third author of this paper came to similar conclusions a couple of years earlier, and told it to the first author, but did not publish the result.)

We slightly improve the results of Lang and Schroeder.
Our proof is based on the barycentric maps introduced by Kleiner in \cite{kleiner}.
The material of this paper will be included in the book on Alexandrov geometry that we are currently writing, but it seems useful to publish it now.

 \parbf{Structure of the paper.}
We introduce notations in Section~\ref{sec:prelim}.
In section \ref{sec:4pt} we give altternative definitions of Alexandrov spaces 
based on the Kirszbraun property for 4-point sets.
The generalized Kirszbraun theorem is proved in Section~\ref{sec:kirszbraun}.
In the sections \ref{sec:1+n} and \ref{sec:2n+2} we describe some comparison properties of finite subsets of Alexandrov spaces.
In Section~\ref{sec:kirszbraun:open} we discuss related open problems.
Appendices~\ref{sec:baricentric} and \ref{sec:helly} describe Kleiner's barycentric map and an analog of Helly's theorem for Alexandrov spaces.

\parbf{Historical remark. }
 Not much is known about the author of this remarkable theorem.
The theorem appears in Kirszbraun's master's thesis which he defended in Warsaw University in 1930.
His name is Moj\.{z}esz and his second name is likely to be Dawid but is uncertain.
He was born either in 1903 or 1904 and died in a ghetto in 1942.
After university, 
he worked  as an actuary in an insurance company; 
\cite{kirszbraun} seems to be his only publication in mathematics. 

\parbf{Acknowledgment.} 
We want to thank S.~Ivanov, N.~Lebedeva and A.~Lytchak 
for useful comments
and pointing out misprints.
Also we want to thank L.~Grabowski for bringing to our attention the entry about Kirszbraun in the Polish Biographical Dictionary.

\section{Preliminaries}\label{sec:prelim}

In this section we mainly introduce our notations.

\parbf{Metric spaces.} Let $\spc{X}$ be a metric space.  The distance between two points $x,y\in\spc{X}$ will be denoted as $\dist{x}{y}{}$ or $\dist{x}{y}{\spc{X}}$.

Given $R\in[0,\infty]$ and $x\in \spc{X}$, the sets
\begin{align*}
\oBall(x,R)&=\{y\in \spc{X}\mid \dist{x}{y}{}<R\},
\\
\cBall[x,R]&=\{y\in \spc{X}\mid \dist{x}{y}{}\le R\}.
\end{align*}
are called respectively the \emph{open} and \emph{closed  ball} of radius $R$ with center at $x$.

A metric space $\spc{X}$ is called 
\emph{intrinsic}
if for any $\eps>0$ and any two points $x,y\in \spc{X}$ with $\dist{x}{y}{}<\infty$ there is an $\eps$-midpoint for $x$ and $y$;
i.e. there is a point $z\in \spc{X}$ such that $\dist{x}{z}{},\dist{z}{y}{}<\tfrac{1}{2}\cdot \dist[{{}}]{x}{y}{}+\eps$.

\parbf{Model space.}
$\Lob{m}{\kappa}$ denotes $m$-dimensional model space with curvature $\kappa$; 
i.e. the simply connected $m$-dimensional Riemannian manifold with constant sectional curvature $\kappa$.

Set $\varpi\kappa=\diam\Lob2\kappa$\index{$\varpi\kappa$}, so 
$\varpi\kappa=\infty$ if $\kappa\le0$ and $\varpi\kappa=\pi/\sqrt{\kappa}$ if $\kappa>0$.
(The letter $\varpi{}$ is a glyph variant of lower case $\pi$,
but is usually pronounced as \emph{pomega}.)

\parbf{Ghost of Euclid.} Let $\spc{X}$ be a metric space 
and $\II$ be a real interval.
A globally isometric map $\gamma\:\II\to \spc{X}$ will be called a \emph{unitspeed geodesic}. 
A unitspeed geodesic between $p$ and $q$ will be denoted by $\geod_{[p q]}$.
We consider $\geod_{[p q]}$ with parametrization starting at $p$; 
i.e. $\geod_{[p q]}(0)=p$ and $\geod_{[p q]}(\dist{p}{q}{})=q$.
The image of $\geod_{[p q]}$ will be denoted by $[p q]$ and called a \emph{geodesic}\index{geodesic}.

Also we will use the following short-cut notation:
\begin{align*}
\l] p q \r[&=[p q]\backslash\{p,q\},
&
\l] p q \r]&=[p q]\backslash\{p\},
&
\l[ p q \r[&=[p q]\backslash\{q\}.
\end{align*}

A metric space $\spc{X}$ is called 
\emph{geodesic}
if for any two points $x,y\in \spc{X}$ there is a geodesic $[x y]$ in $\spc{X}$.

Given a geodesic $[p q]$, we denote by $\dir{p}{q}$ its direction at $p$.
We may think of $\dir{p}{q}$ as belonging to the space of directions $\Sigma_p$ at $p$,
which in turn can be identified with the unit sphere in the tangent space $\T_p$ at $p$.
Further we set $\ddir{p}{q}=\dist[{{}}]{p}{q}{}\cdot\dir{p}{q}$;
it is a \emph{tangent vector} at $p$, that is, an element of $\T_p$.

For a triple of points $p,q,r\in \spc{X}$, a choice of triple of geodesics $([q r], [r p], [p q])$ will be called a \emph{triangle} and we will use the notation 
$\trig p q r=([q r], [r p], [p q])$.
If $p$ is distinct from $x$ and $y$, a  pair of geodesics $([p x],[p y])$ will be called a \emph{hinge}\index{hinge}, and  denoted by 
$\hinge p x y=([p x],[p y])$.

\parbf{Functions.}
A locally Lipschitz function $f$ on a metric space $\spc{X}$ is called $\lambda$-convex ($\lambda$-concave)
if for any geodesic $\geod_{[p q]}$ in $\spc{X}$ the real-to-real function 
$$t\mapsto f\circ\geod_{[p q]}(t)-\tfrac\lambda2\cdot t^2$$
is convex (respectively concave).
In this case we write $f''\ge \lambda$ (respectively $f''\le \lambda$).

A function $f$ is called \emph{strongly convex} (\emph{strongly concave})
if $f''\ge \delta$ (respectively $f''\le -\delta$) for some $\delta>0$.

\parbf{Model angles and triangles.} 
Let $\spc{X}$ be a metric space, 
$p,q,r\in \spc{X}$ 
and $\kappa\in\RR$. 
Let us define a \emph{model triangle} $\trig{\~p}{\~q}{\~r}$ 
(briefly, 
$\trig{\~p}{\~q}{\~r}=\modtrig\kappa(p q r)$) to be a triangle in the model plane $\Lob2\kappa$ such that
$$\dist{\~p}{\~q}{}=\dist{p}{q}{},
\ \ \dist{\~q}{\~r}{}=\dist{q}{r}{},
\ \ \dist{\~r}{\~p}{}=\dist{r}{p}{}.$$
If $\kappa\le 0$, the model triangle is said to be defined, since such a triangle always exists and is unique up to an isometry of $\Lob2\kappa$.
If $\kappa>0$, the model triangle is said to be defined if in addition
$$\dist{p}{q}{}+\dist{q}{r}{}+\dist{r}{p}{}< 2\cdot\varpi\kappa.$$
In this case the triangle also exists and is unique up to an isometry of $\Lob2\kappa$.

If for  $p,q,r\in \spc{X}$, the model triangle 
$\trig{\~p}{\~q}{\~r}=\modtrig\kappa(p q r)$ is defined 
and $\dist{p}{q}{},\dist{p}{r}{}>0$, then the  angle measure of 
$\trig{\~p}{\~q}{\~r}$ at $\~p$ will be called the \emph{model angle} of the triple $p$, $q$, $r$, and will be denoted by
$\angk\kappa p q r$.

\parbf{Curvature bounded below.}
We will denote by $\CBB{}{\kappa}$, complete intrinsic spaces  $\spc{L}$ with curvature $\ge\kappa$ in the sense of Alexandrov.
Specifically,  $\spc{L}\in \CBB{}{\kappa}$ if for any quadruple of points $p,x^1,x^2,x^3\in \spc{U}$ , we have
$$\angk\kappa p{x^1}{x^2}
+\angk\kappa p{x^2}{x^3}
+\angk\kappa p{x^3}{x^1}\le 2\cdot\pi.\eqlbl{Yup-kappa}$$
or at least one of the model angles $\angk\kappa p{x^i}{x^j}$ is not defined.

Condition \ref{Yup-kappa} will be called \emph{(1+3)-point comparison}.

According to Plaut's theorem \cite[Th. 27]{plaut:survey},
any space $\spc{L}\in \CBB{}{}$ is $G_\delta$-geodesic; 
that is, for any point $p\in \spc{L}$ there is a dense $G_\delta$-set $W_p\subset\spc{L}$ such that for any $q\in W_p$ there is a geodesic $[p q]$.

We will use two more equivalent definitions of $\CBB{}{}$ spaces (see \cite{AKP}).
Namely, a complete $G_\delta$-geodesic space is in $\CBB{}{}$ 
if and only if it satisfies either of following conditions: 
\begin{enumerate}

\item\label{POS-CBB-ref} (point-on-side comparison)
For any geodesic $[x y]$ and $z\in \l]x y\r[$, we have
$$\angk\kappa x p y\le\angk\kappa x p z; \eqlbl{POS-CBB}$$
or, equivalently, 
$$\dist{\~p}{\~z}{}\le \dist{p}{z}{},$$
where $\trig{\~p}{\~x}{\~y}=\modtrig\kappa(p x y)$, $\~z\in\l] \~x\~y\r[$, $\dist{\~x}{\~z}{}=\dist{x}{z}{}$.

\item (hinge comparison)
For any hinge $\hinge x p y$, the angle 
$\mangle\hinge x p y$ is defined and 
$$\mangle\hinge x p y\ge\angk\kappa x p y.$$
Moreover, if $z\in\l]x y\r[$, $z\not=p$ then for any two hinges $\hinge z p y$ and $\hinge z p x$ with common side $[z p]$
$$\mangle\hinge z p y + \mangle\hinge z p x\le\pi.$$
\end{enumerate}

We also use the following standard result in Alexandrov geometry, 
which follows from the discussion in the survey of Plaut \cite[8.2]{plaut:survey}.

\begin{thm}{Theorem}\label{thm:cbb-lin-part}
Let $\spc{L}\in \CBB{}{}$.
Given an array of points $(x^1,x^2\dots,x^n)$ in $\spc{L}$,
 there is a dense $G_\delta$-set $W\subset\spc{L}$ such that
for any $p\in W$, all the directions $\dir{p}{x^i}$ lie in 
an isometric copy of a unit sphere in $\Sigma_p$.
(Or, equivaletntly, all the vectors $\ddir{p}{x^i}$ lie in 
a subcone of the tangent space $\T_p$ which is isometric to Euclidean space.)
\end{thm}

\parbf{Curvature bounded above.}
We will denote by $\Cat{}{\kappa}$ 
the class of metric spaces $\spc{U}$ in which any two points at distance $<\varpi\kappa$ are joined by a geodesic, 
and which have  curvature $\le\kappa$ in the following global sense of Alexandrov:  namely, for any quadruple of points $p^1,p^2,x^1,x^2\in \spc{U}$, we have
$$
 \angk{\kappa}{p^1}{x^1}{x^2} 
\le 
\angk{\kappa}{p^1}{p^2}{x^1}+\angk{\kappa}{p^1}{p^2}{x^2},
\ \t{or}\ 
\angk{\kappa} {p^2}{x^1}{x^2}\le \angk{\kappa} {p^2}{p^1}{x^1} + \angk{\kappa} {p^2}{p^1}{x^2},
\eqlbl{gokova:eq:2+2}$$
 or
one of the six model angles above 
is undefined.

The condition \ref{gokova:eq:2+2} will be called \emph{(2+2)-point comparison} (or \emph{(2+2)-point $\kappa$-comparison} 
if a confusion may arise).  

We denote the complete $\Cat{}{\kappa}$ spaces by $\cCat{}{\kappa}$. 

The following lemma is a direct consequence of the definition:

\begin{thm}{Lemma}\label{lem:cat-complete} 
Any complete intrinsic 
 space $\spc{U}$  in which every quadruple $p^1,p^2,x^1,x^2$ satisfies 
the (2+2)-point $\kappa$-
comparison 
is   a $\cCat{}{\kappa}$ space (that is, any  two points at distance $<\varpi\kappa$ are  joined by a geodesic).  

In particular,  the completion of a 
$\Cat{}{\kappa}$ space again lies in $\Cat{}{\kappa}$.
\end{thm}


We have the following basic facts (see [1]):

\begin{thm}{Lemma}\label{lem:cat-unique} 
In a $\Cat{}{\kappa}$ space, geodesics of length $<\varpi\kappa$ are uniquely determined by, and continuously dependent on, their endpoint pairs.\end{thm}

\begin{thm}{Lemma}\label{lem:convex-balls}
In a $\Cat{}{\kappa}$ space, any open ball  $\oBall(x,R)$ of radius $R\le\varpi\kappa/2$  is convex, that is, $\oBall(x,R)$ contains every geodesic whose endpoints it contains.
\end{thm}

We also use an equivalent definition of $\Cat{}{\kappa}$ spaces (see \cite{AKP}).
Namely, a metric space  $\spc{U}$ in which any  two points at distance $<\varpi\kappa$ are  joined by a geodesic is a $\Cat{}{\kappa}$ space if and only if it satisfies the following condition: 
\begin{enumerate}

\item (point-on-side comparison)\label{cat-monoton}
for any geodesic $[x y]$ and $z\in \l]x y\r[$, we have
$$\angk\kappa x p y\ge\angk\kappa x p z,$$
or equivalently, 
$$\dist{\~p}{\~z}{}\ge \dist{p}{z}{}, \eqlbl{POS-CAT}$$
where $\trig{\~p}{\~x}{\~y}=\modtrig\kappa(p x y)$, $\~z\in\l] \~x\~y\r[$, $\dist{\~x}{\~z}{}=\dist{x}{z}{}$.


\end{enumerate}

We also use Reshetnyak's majorization theorem \cite{reshetnyak:major}.  
Suppose $\~\alpha$ is a simple closed curve of finite length  in $\Lob2{\kappa}$,
and $D\subset\Lob2{\kappa}$ is a closed region bounded by $\~\alpha$. If $\spc{X}$ is a metric space,  a length-nonincreasing map $F\:D\to\spc{X}$ is called \emph{majorizing} if it is length-preserving on $\~\alpha$.
In this case, we say that $D$ \emph{majorizes} the curve $\alpha=F\circ\~\alpha$ under the map $F$.

\begin{thm}{Reshetnyak's majorization theorem}
\label{thm:major}
Any closed curve $\alpha$ of length $<2\cdot \varpi\kappa$ in $\spc{U}\in\Cat{}{\kappa}$ is majorized by a convex region in $\Lob2\kappa$.
\end{thm}

\parbf{Ultralimit of metric spaces.}
Given a metric space $\spc{X}$, its ultrapower 
(i.e. ultralimit of constant sequence $\spc{X}_n=\spc{X}$) will be denoted as $\spc{X}^\o$;
here $\o$ denotes a fixed nonprinciple ultrafilter.
For definitions and properties of ultrapowers, 
we refer to a paper of Kleiner and Leeb \cite[2.4]{kleiner-leeb}.

We use the following facts about ultrapowers which easily follow from the definitions (see \cite{AKP} for details):
\begin{itemize}
\item $\spc{X}\in\cCat{}{\kappa}\  \Longleftrightarrow\  \spc{X}^\o\in\cCat{}{\kappa}$.
\item $\spc{X}\in\CBB{}{\kappa}\  \Longleftrightarrow\  \spc{X}^\o\in\CBB{}{\kappa}$.
\item $\spc{X}$ is intrinsic if and only if $\spc{X}^\o$ is geodesic.
\end{itemize}

Note that if $\spc{X}$ is \textit{proper} (namely, bounded closed sets are compact), then $\spc{X}$ and $\spc{X}^\o$ coincide.
Thus a reader interested only in proper spaces may ignore everything related to ultrapower in this article.

\section{Short map extension definitions.}\label{sec:4pt}

Theorems \ref{thm:kirsz-def} and \ref{thm:cba-kirsz-def} 
 give characterizations of $\CBB{}{\kappa}$ and  $\Cat{}{\kappa}$.
Very similar theorems were proved by Lang and Shroeder in 
\cite{lang-schroeder}.

\begin{thm}{Theorem}\label{thm:kirsz-def} 
Let $\spc{L}$ be a complete intrinsic space. 
Then $\spc{L}\in\CBB{}{\kappa}$ if and only if for any 3-point set $V_3$ and any 4-point set $V_4\supset V_3$ in $\spc{L}$, 
any short map $f\:V_3\to\Lob2\kappa$ can be extended to a short map $F\:V_4\to\Lob2\kappa$ (so $f=F|_{V_3}$).
\end{thm}

\begin{thm}{Theorem}\label{thm:cba-kirsz-def} 
Let $\spc{U}$ be a metric space in which any pair of points at distance $<\varpi\kappa$ are joined by a unique geodesic. Then $\spc{U}\in\Cat{}{\kappa}$ if and only if for any $3$-point set $V_3$ and  $4$-point set $V_4\supset V_3$ in $\Lob2\kappa$, where the perimeter of $V_3$ is $<2\cdot\varpi\kappa$, any short map $f\:V_3\to\spc{U}$ can be extended to a short map $F\:V_4\to\spc{U}$.
\end{thm}

The proof of the ``only if'' part of Theorem \ref{thm:kirsz-def} can be obtained as a corollary of Kirszbraun's theorem (\ref{thm:kirsz+}).
But we present another proof, based on more elementary ideas.  The ``only if'' part of Theorem \ref{thm:cba-kirsz-def} does not follow directly from Kirszbraun's theorem, since the desired extension is in $\spc{U}$, not just the completion of $\spc{U}$.

In the proof of Theorem~\ref{thm:kirsz-def}, we use the following lemma in the  geometry of model planes.
Here we say that  two triangles with a common vertex  \emph{do not overlap} if their convex hulls intersect only at the common vertex.

\begin{thm}{Overlap lemma} \label{lem:extend-overlap}
Let $\trig{\~x^1}{\~x^2}{\~x^3}$ be a triangle in $\Lob2{\kappa}$.  Let $\~p^1,\~p^2,\~p^3$ be points such that, for any permutation $\{i,j,k\}$ or $\{1,2,3\}$, we have
\begin{enumerate}[(i)]

\item 
\label{no-overlap:px=px}
$\dist{\~p^i}{\~x^\kay}{}=\dist{\~p^j}{\~x^\kay}{}$,

\item
\label{no-overlap:orient-1}
$\~p^i$ and $\~x^i$ lie in the same closed halfspace determined by $[\~x^j\~x^\kay]$,  

\item
\label{no-overlap:orient-2}
$\mangle\hinge {\~x^i}{\~x^j}{\~p^k}+\mangle\hinge {\~x^i}{\~p^j}{\~x^k}<\pi$.
\end{enumerate}
Set $\mangle\~p^i=\mangle\hinge{\~p^i}{\~x^\kay}{\~x^j}$. It follows that:
\begin{subthm}{two-overlap}
If $\mangle{\~p^1} +\mangle {\~p^2}+\mangle {\~p^3} \le 2\cdot\pi$ and 
triangles $\trig{\~p^3}{\~x^1}{\~x^2}$, $\trig{\~p^2}{\~x^3}{\~x^1}$ do not overlap, then
$$ \mangle {\~p^1} > \mangle{\~p^2}+ \mangle{\~p^3}.$$
\end{subthm}
\begin{subthm}{no-overlap} 
No pair of triangles $\trig{\~p^i}{\~x^j}{\~x^\kay}$  overlap if and only if 
$$\mangle{\~p^1} +\mangle {\~p^2}+\mangle{\~p^3}> 2\cdot\pi.$$
\end{subthm}
\end{thm}

\begin{wrapfigure}{r}{20mm}
\begin{lpic}[t(-0mm),b(0mm),r(0mm),l(0mm)]{pics/contr-no-overlap(0.4)}
\lbl[rt]{37,34;$\~p^1$}
\lbl[tl]{10,36;$\~p^2$}
\lbl[bl]{14,6;$\~p^3$}
\lbl[lb]{13,51;$\~x^1$}
\lbl[tr]{12,0;$\~x^2$}
\lbl[l]{44,30;$\~x^3$}
\end{lpic}
\end{wrapfigure}

\parbf{Remark.}
If $\kappa\le 0$, the ``only if'' part of (\ref{SHORT.no-overlap}) can be proved without using condition (\ref{no-overlap:px=px}).
This follows immediately from the formula that relates the sum of angles for the hexagon
$[\~p^1\~x^2\~p^3\~x^1\~p^2\~x^3]$ and its area:
$$ \mangle\~p^1
-
\mangle\~x^2
+
\mangle\~p^3
-
\mangle\~x^1
+
\mangle\~p^2
-
\mangle\~x^3
=2\cdot\pi-\kappa\cdot{\area}.
$$

In case $\kappa>0$, condition (\ref{no-overlap:px=px}) is essential.
An example for $\kappa>0$ can be constructed by perturbing the degenerate spherical configuration on the picture.

\parit{Proof.}   Rotate the  triangle $\trig{\~p^3}{\~x^1}{\~x^2}$ around $\~x^1$ to make $[\~x^1\~p^3]$ coincide with $[\~x^1\~p^2]$.
Let  $\dot x^2$ denote the image of $\~x^2$ after rotation. 
By (\ref{no-overlap:orient-1}) and (\ref{no-overlap:orient-2}), 
the triangles $\trig{\~p^3}{\~x^1}{\~x^2}$ and $\trig{\~p^2}{\~x^3}{\~x^1}$ do not overlap if and only if 
$\mangle\hinge{\~x^1}{\~x^3}{\dot x^2}
< \mangle\hinge{\~x^1}{\~x^3}{\~x^2}$, and hence if and only if $\dist{\dot x^2}{\~x^3}{} <   \dist{\~x^2}{\~x^3}{}$. This inequality holds if and only if 
$$
\begin{aligned}
\mangle\~p^1
&> \mangle\hinge{\~p^2}{\~x^3}{\dot x^2}
\\
&=
\min\{\mangle\~p^3+\mangle\~p^2,2\cdot\pi -(\mangle\~p^3+\mangle\~p^2)\},
\end{aligned}
\eqlbl{eq:no-overlap}$$
since  in the  inequality, the  corresponding hinges have the same pairs of sidelengths.
(The two pictures show that both possibilities for the minimum can occur.)

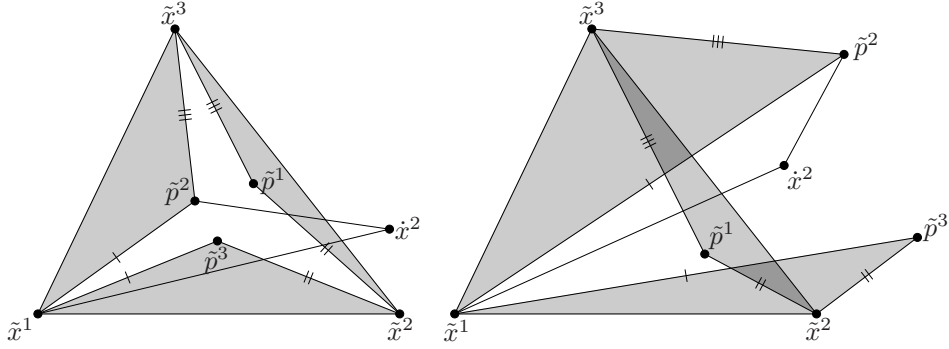
\begin{center} 
\begin{lpic}[t(0mm),b(10mm),r(0mm),l(0mm)]{pics/4-pnt-kirsz-x2(0.2)}
\lbl[t]{125,52;$\~p^3$}
\lbl[l]{154,96;$\~p^1$}
\lbl[rb]{108,80;$\~p^2$}
\lbl[tr]{5,5;$\~x^1$}
\lbl[t]{246,5;$\~x^2$}
\lbl[l]{242,65;$\dot x^2$}
\lbl[b]{97,200;$\~x^3$}
\lbl[lb]{449,51;$\~p^1$}
\lbl[l]{544,185;$\~p^2$}
\lbl[l]{590,65;$\~p^3$}
\lbl[t]{283,4;$\~x^1$}
\lbl[t]{521,5;$\~x^2$}
\lbl[lt]{500,105;$\dot x^2$}
\lbl[b]{372,200;$\~x^3$}
\end{lpic}

\end{center}

If $\mangle\~p^1 + \mangle\~p^2+\mangle\~p^3 \le 2\cdot\pi$,
then  \ref{eq:no-overlap} implies $\mangle\~p^1>\mangle\~p^2 + \mangle\~p^3$.
That proves (\ref{SHORT.two-overlap}).

\parit{``Only if'' part of (\ref{SHORT.no-overlap}).} 
Suppose no two triangles overlap and $\mangle\~p^2 + \mangle\~p^2+\mangle\~p^3 \le 2\cdot\pi$.  By \ref{SHORT.two-overlap}), for $\{i,j,k\}=\{1,2,3\}$ we have
$$\mangle\~p^i > \mangle\~p^j+\mangle\~p^\kay.$$ 
Adding these three inequalities gives a contradiction:
$$ \mangle\~p^1+\mangle\~p^2+\mangle\~p^3 > 2\cdot (\mangle\~p^1+\mangle\~p^2+\mangle\~p^3).$$

\parit{``If'' part of (\ref{SHORT.no-overlap}). } 
Suppose triangles $\trig{\~p^3}{\~x^1}{\~x^2}$ and $\trig{\~p^2}{\~x^3}{\~x^1}$ overlap 
and 
$$\mangle\~p^1 + \mangle\~p^2+\mangle\~p^3 > 2\cdot\pi.
\eqlbl{eq:<p1+<p2+<p3}$$  
By the former, \ref{eq:no-overlap} fails.  
By \ref{eq:<p1+<p2+<p3}, $\mangle\~p^2+\mangle\~p^3 > \pi$. 
Therefore
$$\mangle\~p^1
\le 2\cdot\pi -(\mangle\~p^2+\mangle\~p^3),$$ 
which contradicts \ref{eq:<p1+<p2+<p3}.
\qeds

\parit{Proof of \ref{thm:kirsz-def}; ``if'' part.} 
Assume $\spc{L}$ is geodesic.
Let  $x^1,x^2,x^3\in \spc{L}$ be such that the model triangle 
$\trig{\~x^1}{\~x^2}{\~x^3}=\modtrig\kappa(x^1 x^2 x^3)$ is defined.
Choose $p\in \,{]}x^1x^2{[}\,$.
Let  $V_3=\{x^1,x^2,x^3\}$ and 
$V_4=\{x^1,x^2,x^3,p\}$, and set  $f(x^i)=\~x^i$.  Then a short extension of $f$ to $V_4$ gives point-on-side comparison (see page~\pageref{POS-CBB}).

In case $\spc{L}$ is not geodesic, pass to its ultrapower $\spc{L}^\o$.
Note that if $\spc{L}$ satisfies the conditions of Theorem \ref{thm:kirsz-def}  then so does $\spc{L}^\o$. Also, 
recall that $\spc{L}^\o$ is geodesic. 
Thus, from above, ${\spc{L}^\o}\in\CBB{}{\kappa}$.
Hence $\spc{L}\in\CBB{}{\kappa}$.

\parit{``Only if'' part.}
Assume the contrary;
i.e.,  $x^1,x^2,x^3,p\in \spc{L}$, and 
$\~x^1,\~x^2,\~x^3\in\Lob2\kappa$ are such that
$\dist{\~x^i}{\~x^j}{}\le\dist{x^i}{x^j}{}$ for all $i,j$ and there is no point $\~p\in \Lob2\kappa$ such that $\dist{\~p}{\~x^i}{}\le \dist{p}{x^i}{}$ for all $i$.

We claim that in this case all comparison triangles $\modtrig\kappa(p x^ix^j)$ are defined.
That is always true if $\kappa\le0$.
If $\kappa>0$, and say $\modtrig\kappa(p x^1x^2)$ is undefined, then 
\begin{align*}
\dist{p}{x^1}{}+\dist{p}{x^2}{}
&\ge 
2\cdot\varpi\kappa-\dist{x^1}{x^2}{}
\ge
\\
&\ge 
2\cdot\varpi\kappa-\dist{\~x^1}{\~x^2}{}
\ge
\\
&\ge 
\dist{\~x^1}{\~x^3}{}+\dist{\~x^2}{\~x^3}{}.
\end{align*}
Thus one can take $\~p$ on $[\~x^1\~x^3]$ or  $[\~x^2\~x^3]$.

For each $i\in \{1,2,3\}$, consider a point $\~p^i\in\Lob2\kappa$ such that $\dist{\~p^i}{\~x^i}{}$ is minimal among points satisfying $\dist{\~p^i}{\~x^j}{}\le\dist{p}{ x^j}{}$ for all $j\not=i$. 
Clearly, every $\~p^i$ is inside the triangle $\trig{\~x^1}{\~x^2}{\~x^3}$ (that is, in $\Conv(\~x^1,\~x^2,\~x^3)$), and $\dist{\~p^i}{\~x^i}{}>\dist{p}{ x^i}{}$ for each $i$.
It follows that
\begin{enumerate}[(i)]
\item $\dist{\~p^i}{\~x^j}{}=\dist{p}{ x^j}{}$ for $i\not=j$;
\item no pair of triangles from $\trig{\~p^1}{\~x^2}{\~x^3}$, $\trig{\~p^2}{\~x^3}{\~x^1}$, $\trig{\~p^3}{\~x^1}{\~x^2}$ overlap in $\trig{\~x^1}{\~x^2}{\~x^3}$.
\end{enumerate}

As follows from Lemma~\ref{no-overlap}, in this case 
$$\mangle\hinge {\~p^1}{\~x^2}{\~x^3} 
+\mangle\hinge {\~p^2}{\~x^3}{\~x^1}
+\mangle\hinge {\~p^3}{\~x^1}{\~x^2}
>2\cdot\pi.$$
Thus we arrive at a contradiction, since $\dist{\~x^i}{\~x^j}{}\le\dist{x^i}{x^j}{}$ implies that
$$\mangle\hinge {\~p^\kay}{\~x^i}{\~x^j}
\le
\angk\kappa p{x^i}{x^j}$$
if $(i,j,k)$ is a permutation of $(1,2,3)$.
\qeds

In the proof of Theorem~\ref{thm:cba-kirsz-def},
we use the following lemma in the geometry of model planes: 

\begin{thm}{Lemma}\label{lem:smaller-trig}
Let $x^1,x^2,x^3,y^1,y^2,y^3\in\Lob{}{\kappa}$
be points such that $\dist{x^i}{x^j}{}\ge\dist{y^i}{y^j}{}$ for all $i,j$.
Then there is a short map $\map\:\Lob{}{\kappa}\to\Lob{}{\kappa}$ such that $\map(x^i)=y^i$ for all $i$;
moreover, one can choose $\map$ so that 
$$\Im \map\subset\Conv(y^1,y^2,y^3).$$

\end{thm}

We only give an idea of the proof of  this lemma;
alternatively, one can get the result as a corollary of  Kirszbraun's theorem (\ref{thm:kirsz+}) 

\parit{Idea of the proof.}
The map $\map$ can be constructed as a composition of the following folding maps:
Given a halfspace $H$ in $\Lob{}{\kappa}$, consider the map $\Lob{}{\kappa}\to H$, 
which is the identity on $H$ and reflects all points outside of $H$ into $H$.
This map is a path isometry, in particular, it is short. 

One can get the last part of the lemma by composing the above map with foldings along the sides of triangle $\trig{y^1}{y^2}{y^3}$ and passing to a partial limit.
\qeds

\parit{Proof of \ref{thm:cba-kirsz-def}; ``if'' part.}
The point-on-side comparison (\ref{cat-monoton}) follows  by
taking $V_3=\{\~x,\~y,\~p\}$ and  $V_4=\{\~x,\~y,\~p,\~z\}$ where $z\in \l]x y\r[$.  
It is only necessary to observe that  $F(\~z)=z$ by uniqueness of $[x y]$.

\parit{``Only if'' part.}
Let $V_3=\{\~x^1,\~x^2,\~x^3\}$ and $V_4=\{\~x^1,\~x^2,\~x^3,\~p\}$.

Set $y^i\z=f(\~x^i)$ for all $i$;
we need to find a point $q\in\spc{U}$ such that $\dist{y^i}{q}{}\le\dist{\~x^i}{\~p}{}$ for all $i$.

Consider the model triangle $\trig{\~y^1}{\~y^2}{\~y^3}=\modtrig\kappa({y^1}{y^2}{y^3})$. 
Set $D\z=\Conv({\~y^1},{\~y^2},{\~y^3})$.

Note that $\dist{\~y^i}{\~y^j}{}=\dist{y^i}{y^j}{}\le\dist{\~x^i}{\~x^j}{}$ for all $i,j$.
Applying Lemma \ref{lem:smaller-trig},
we get a short map 
$\map\:\Lob{}{\kappa}\to D$ such that 
$\map\:\~x^i\mapsto\~y^i$.

Further, from Reshetnyak majorization (\ref{thm:major}), 
there is a short map $F\:D\to \spc{U}$ such that $\~y^i\mapsto y^i$ for all $i$.

Thus one can take $q=F\circ\map(\~p)$.
\qeds

\section{(1+\textit{n})-point comparison}\label{sec:1+n}

The following theorem gives a more sensitive analog of (1+3)-point comparison.
In a bit more analytic form it was discovered by Sturm in \cite{sturm}.

\begin{thm}{(1+\textit{n})-point comparison}
\label{thm:pos-config} 
Let $\spc{L}\in\CBB{}{\kappa}$.
Then for any array of points $p,x^1,\dots,x^n\in \spc{L}$  
there is a model array $\~p,\~x^1,\dots,\~x^n\in\Lob{n}\kappa$ such that
\begin{subthm}{}
$\dist{\~p}{\~x^i}{}=\dist{p}{x^i}{}$ for all $i$.
\end{subthm}

\begin{subthm}{}$\dist{\~x^i}{\~x^j}{}\ge\dist{x^i}{x^j}{}$ for all $i,j$.
\end{subthm}
\end{thm}

\parit{Proof.} 
It is enough to show that given $\eps>0$ there is a configuration $\~p,\~x^1,\dots,\~x^n\in\Lob{n}\kappa$ such that $\dist{\~x^i}{\~x^j}{}\ge\dist{x^i}{x^j}{}$ and $\bigl|\dist{\~p}{\~x^i}{}-\dist{p}{x^i}{}\bigr|\le \eps$.
Then one can pass to a limit configuration for $\eps\to 0+$.

According to \ref{thm:cbb-lin-part}, there is a point $p'$ such  that $\dist{p'}{p}{}\le\eps$
and
$\T_{p'}$ contains a subcone $E$ isometric to a Euclidean space 
which contains all vectors $\ddir{p'}{x^i}$.
Passing to a subspace if necessary, we can assume that $\dim E\le n$.

Mark a point $\~p\in \Lob{n}\kappa$ and choose an isometric embedding
$\imath\: E\to \T_{\~p}\Lob{n}\kappa$.
Set 
$$\~x^i=\exp_{\~p}\circ\imath\circ\ddir{p'}{x^i}.$$
Thus $\dist{\~p}{\~x^i}{}=\dist{p'}{x^i}{}$ and therefore
$$\bigl|\dist{\~p}{\~x^i}{}-\dist{p}{x^i}{}\bigr|\le \dist{p}{p'}{} \le\eps.$$
From the hinge comparison,
we have 
$$\angk\kappa{\~p}{\~x^i}{\~x^j}
=\mangle\hinge{\~p}{\~x^i}{\~x^j}
=\mangle\hinge{p'}{x^i}{x^j}\ge \angk\kappa{p'}{x^i}{x^j},$$
thus 
$$\dist{\~x^i}{\~x^j}{}\ge \dist{x^i}{x^j}{}.$$
\qedsf

\section{Kirszbraun's theorem}\label{sec:kirszbraun}

A slightly weaker version of the following theorem
was proved by Lang and Schroeder in \cite{lang-schroeder}.
The Conjecture~\ref{conj:kirsz} (if true) gives an equivalent condition for the existence of a short extension;
roughly it states that example \ref{example:SS_+} is the only obstacle.

\begin{thm}{Kirszbraun's theorem}
\label{thm:kirsz+}
Let
$\spc{L}\in\CBB{}{\kappa}$, 
$\spc{U}\in\cCat{}{\kappa}$, 
$Q\subset \spc{L}$ be arbitrary subset
and $f\: Q\to\spc{U}$ be a short map.
Assume that there is $z\in\spc{U}$ such that 
$f(Q)\subset \oBall[z,\tfrac{\varpi\kappa}{2}]$.
Then $f\:Q\to\spc{U}$ can be extended to a short map 
$F\:\spc{L}\to \spc{U}$
(that is, there is a short map $F\:\spc{L}\to \spc{U}$ such that $F|_Q=f$.)
\end{thm}
 
The condition $f(Q)\subset \oBall[z,\tfrac{\varpi\kappa}{2}]$ trivially holds for any $\kappa\le 0$ since in this case $\varpi\kappa=\infty$. 
The following example shows that this condition is needed for $\kappa>0$.

\begin{thm}{Example}\label{example:SS_+}
Let $\SS^m_+$ be a closed $m$-dimensional unit hemisphere.  Denote its boundary, which is isometric to $\SS^{m-1}$, by  $\partial\SS^m_+$.
Clearly, $\SS^m_+\in\CBB{}1$ and $\partial\SS^m_+\in\cCat{}1$ but the identity map ${\partial\SS^m_+}\to \partial\SS^m_+$ cannot be extended to a short map $\SS^m_+\to \partial\SS^m_+$ (there is no place for the pole).

There is also a direct generalization of this example to a hemisphere in a Hilbert space of arbitrary cardinal dimension.
\end{thm}

First we prove this theorem in the case $\kappa\le 0$ (\ref{thm:kirsz}).
In the proof of the more complicated case $\kappa>0$, we use the case $\kappa=0$.
The following lemma is the main ingredient in the proof. 

\begin{thm}{Finite$\bm{+}$one lemma}\label{lem:kirsz-neg:new}
Let $\kappa\le 0$,
$\spc{L}\in\CBB{}{\kappa}$, and 
$\spc{U}\in\cCat{}{\kappa}$.  
Let  
$x^1,x^2,\dots,x^n\in\spc{L}$ 
and $y^1,y^2,\dots,y^n\in\spc{U}$
be
such that $\dist{x^i}{x^j}{}\ge\dist{y^i}{y^j}{}$ for all $i,j$.

Then for any $p\in\spc{L}$, there is $q\in\spc{U}$ such that $\dist{y^i}{q}{}\le\dist{x^i}{p}{}$ for each $i$.
\end{thm}

\parit{Proof.}
It is sufficient to prove the lemma only for $\kappa=0$ and $-1$.
The proofs of these two cases are identical, only the formulas differ.
In the proof, we assume $\kappa=0$ and provide the formulas for $\kappa=-1$ in the footnotes.

From (1+\textit{n})-point comparison (\ref{thm:pos-config}), 
there is a model configuration 
$\~p,\~x^1,\~x^2,\dots,\~x^n\in \Lob{n}{\kappa}$ such that
$\dist{\~p}{\~x^i}{}=\dist{p}{x^i}{}$
and $\dist{\~x^i}{\~x^j}{}\ge\dist{x^i}{x^j}{}$ 
for all $i$, $j$.

For each $i$, consider functions 
$f^i\:\spc{U}\to\RR$ and $\~f^i\:\Lob{n}{\kappa}\to\RR$ 
defined as follows%
\footnote{In case $\kappa=-1$,
$$
\begin{aligned}
&f^i(y)=\cosh\dist{y^i}{y}{},
&
&\~f^i(\~x)=\cosh\dist{\~x^i}{\~x}{}.
\end{aligned}\eqno{(A)\mc-}$$}:
$$
\begin{aligned}
&f^i(y)=\tfrac{1}{2}\cdot\dist[2]{y^i}{y}{},
&
&\~f^i(\~x)=\tfrac{1}{2}\cdot\dist[2]{\~x^i}{\~x}{}.
\end{aligned}\eqno{(A)\mc0}
$$
Set
$\bm{f}=(f^1,f^2,\dots,f^n)\:\spc{U}\to\RR^n$ and $\bm{\~f}=(\~f^1,\~f^2,\dots,\~f^n)\:\Lob{n}{\kappa}\to\RR^n$.

Recall that $\SupSet$ (superset in $\RR^n$) is defined in \ref{def:supset+succcurlyeq}.
Note that it is sufficient to prove that
$\bm{\~f}(\~p)\in\SupSet\bm{f}(\spc{U})$.

Clearly,
$(f^i)''\ge 1$.
Thus, by the theorem on barycentric simplex (\ref{bary-iff}), 
the set $\SupSet\bm{f}(\spc{U})\subset\RR^{n}$ is convex.

Arguing by contradiction, let us assume that $\bm{\~f}(\~p)\not\in\SupSet\bm{f}(\spc{U})$.

Then there  exists a supporting hyperplane  $\alpha_1x_1+\ldots \alpha_nx_n=c$ to $\SupSet\bm{f}(\spc{U})$, separating it from  $\bm{\~f}(\~p)$.
Just as in the proof of Theorem~\ref{thm:bary} we have that all $\alpha_i\ge 0$. So by rescaling we can assume that $(\alpha_1,\alpha_2,\dots,\alpha_n)\in\Delta^{n-1}$ and 
$$\sum_i\alpha_i\cdot\~f^i(\~p)
< 
\inf
\set{\sum_i\alpha_i\cdot f^i(q)}{q\in\spc{U}}.$$ 
The latter contradicts the following claim.

\begin{clm}{Claim}
Given $\bm{\alpha}=(\alpha_1,\alpha_2,\dots,\alpha_n)\in\Delta^{n-1}$,
set
\begin{align*}
&h=\sum_i\alpha_i\cdot f^i,
&
&h\:\spc{U}\to\RR,
&
&z=\argmin h\in \spc{U},
\\
&\~h=\sum_i\alpha_i\cdot \~f^i,
&
&\~h\:\Lob{n}{\kappa}\to\RR,
&
&\~z=\argmin \~h\in \Lob{n}{\kappa}.
\end{align*}
Then 
$h(z)\le \~h(\~z)$.
\end{clm}

\parit{Proof of the claim.}
Note that $\d_z h\ge 0$.
Thus, for each $i$, we have%
\footnote{In case $\kappa=-1$, the same calculations give
$$
\begin{aligned}
0
&\le\dots \le
-\tfrac{1}{\sinh\dist[{{}}]{z}{y^i}{}}
\cdot 
\sum_j
\alpha_j\cdot\l[\cosh\dist[{{}}]{z}{y^i}{}\cdot\cosh\dist[{{}}]{z}{y^j}{}-\cosh\dist[{{}}]{y^i}{y^j}{}\r].
\end{aligned}
\eqno{(B)\mc-}
$$
}
$$
\begin{aligned}
0
&\le (\d_z h)(\dir{z}{y^i})
=
\\
&=
-\sum_j\alpha_j\cdot\dist[{{}}]{z}{y^j}{}\cdot\cos\mangle\hinge{z}{y^i}{y^j}
\le
\\
&\le
-\sum_j\alpha_j\cdot\dist[{{}}]{z}{y^j}{}\cdot\cos\angk0{z}{y^i}{y^j}
=
\\
&=
-\tfrac{1}{2\cdot\dist[{{}}]{z}{y^i}{}}
\cdot 
\sum_j
\alpha_j\cdot\l[\dist[2]{z}{y^i}{}+\dist[2]{z}{y^j}{}-\dist[2]{y^i}{y^j}{}\r].
\end{aligned}
\eqno{(B)\mc0}$$
In particular%
\footnote{In case $\kappa=-1$, the same calculations give
$$
\begin{aligned} 
\sum_{i}\alpha_i\cdot\l[\sum_j
\alpha_j\cdot\l[\cosh\dist[{{}}]{z}{y^i}{}\cdot\cosh\dist[{{}}]{z}{y^j}{}
-\cosh\dist[{{}}]{y^i}{y^j}{}\r]
\r]\le0
\end{aligned}.
\eqno{(C)\mc-}
$$
},
$$
\begin{aligned}
\sum_{i}
\alpha_i
\cdot
\l[\sum_j
\alpha_j
\cdot
\l[\dist[2]{z}{y^i}{}+\dist[2]{z}{y^j}{}-\dist[2]{y^i}{y^j}{}\r]
\r]\le 0,
\end{aligned}
\eqno{(C)\mc0}
$$
or%
\footnote{In case $\kappa=-1$,
$$(h(z))^2\le
\sum_{i,j}
\alpha_i\cdot\alpha_j
\cdot
\cosh\dist[{{}}]{y^i}{y^j}{}. \eqno{(D)\mc-}$$
}
$$2\cdot h(z)
\le
\sum_{i,j}
\alpha_i\cdot\alpha_j
\cdot
\dist[2]{y^i}{y^j}. \eqno{(D)\mc0}$$

Note that if $\spc{U}\iso\Lob{n}{\kappa}$, 
then all inequalities in $(B,C,D)$ are sharp.
Thus the same argument as above, repeated for $\~x^1,\~x^2,\dots,\~x^n\in\Lob{n}{\kappa}$
gives%
\footnote%
{In case $\kappa=-1$,
$$(\~h(\~z))^2
=
\sum_{i,j}
\alpha_i\cdot\alpha_j
\cdot
\cosh\dist[{{}}]{\~x^i}{\~x^j}{}.
\eqno{(E)\mc-}$$
}
$$2\cdot \~h(\~z)
=
\sum_{i,j}
\alpha_i\cdot\alpha_j
\cdot
\dist[2]{\~x^i}{\~x^j}{}. 
\eqno{(E)\mc0}$$
Note that 
$$\dist{\~x^i}{\~x^j}{}\ge\dist{x^i}{x^j}{}\ge\dist{y^i}{y^j}{}$$ 
for all $i$, $j$.
Thus, $(D)$ and $(E)$ imply the claim.
\qeds

\begin{thm}{Kirszbraun's theorem for nonpositive bound}
\label{thm:kirsz}
Let
$\kappa\le0$,
$\spc{L}\in\CBB{}{\kappa}$, 
$\spc{U}\in\cCat{}{\kappa}$, 
$Q\subset \spc{L}$ be arbitrary subset
and $f\: Q\to\spc{U}$ be a short map.
Then there is a short extension 
$F\:\spc{L}\to \spc{U}$ of $f$;
that is, there is a short map $F\:\spc{L}\to \spc{U}$ such that $F|_Q=f$.
\end{thm}

\parbf{Remark.}
If $\spc{U}$ is proper, then in the following proof the Helly's theorem (\ref{thm:helly}) is not needed.
Everything follows directly from compactness of closed balls in $\spc{U}$.

\parit{Proof of \ref{thm:kirsz}.} 
By Zorn's lemma, we can assume 
that $Q\subset\spc{L}$ is a maximal set;
i.e. $f\:Q\to\spc{U}$ does not admits a short extension to any larger set $Q'\supset Q$.

Let us argue by contradiction.
Assume that $Q\not=\spc{L}$;
choose $p\in \spc{L}\backslash Q$.
Then
$$\bigcap_{x\in Q} \cBall[f(x),\dist{p}{x}{}]
=
\emptyset.$$

Since $\kappa\le 0$, the balls are convex; 
thus, by Helly's theorem (\ref{thm:helly}), 
one can choose a point array $x^1,x^2,\dots, x^n\in Q$ such that
$$\bigcap_{i=1}^n \cBall[y^i,\dist{x^i}{p}{}]
=
\emptyset,
\eqlbl{eq:cap=cBalls=0}$$
where $y^i=f(x^i)$.
Finally note that \ref{eq:cap=cBalls=0} contradicts the Finite+one lemma (\ref{lem:kirsz-neg:new}).
\qeds

\parit{Proof of Kirszbraun's theorem (\ref{thm:kirsz+}).} 
The case $\kappa\le 0$ is already proved in \ref{thm:kirsz}.
Thus it remains to prove the theorem only in case $\kappa>0$.
After rescaling we may assume that $\kappa=1$
and therefore $\varpi\kappa=\pi$.

Since $\cBall[z,\pi/2]\in\cCat{}{\kappa}$ (\ref{lem:convex-balls}, \ref{lem:cat-complete}),
we can assume $\spc{U}=\cBall[z,\pi/2]$. 
In particular, any  two points of $\spc{U}$  at distance $<\pi$ are  joined by a geodesic, and $\diam\spc{U}\le\pi$.
If $\dist{x}{y}{}=\pi$ for some $x,y\in\spc{U}$, then the concatenation  of 
$[x z]$ and $[z y]$ forms a geodesic $[x y]$.
Hence $\spc{U}$ is geodesic.

Further, we can also assume that $\diam\spc{L}\le\pi$.
Otherwise $\spc{L}$ is one-dimensional;
in this case the result follows since $\spc{U}$ is geodesic.

\medskip

Assume the theorem is false. Then 
there is a set $Q\subset \spc{L}$, 
a short map $f\: Q\to \spc{U}$ and  
$p\in \spc{L}\backslash  Q$ such that 
$$\bigcap_{x\in  Q}
\cBall[f(x),\dist{x}{p}{}]=\emptyset.
\eqlbl{eq:cap-of-balls}$$

We are going to apply \ref{thm:kirsz} for $\kappa=0$ to the Euclidean cones $\mathring{\spc{L}}=\Cone \spc{L}$ and $\mathring{\spc{U}}\z=\Cone \spc{U}$. 
Note that 
\begin{itemize}
\item ${\mathring{\spc{U}}}\in\cCat{}0$, 
\item since $\diam \spc{L}\le \pi$ we have ${\mathring{\spc{L}}}\in\CBB{}0$. 
\end{itemize}
Further, we view the spaces $\spc{L}$ and $\spc{U}$ as unit spheres in $\mathring{\spc{L}}$ and $\mathring{\spc{U}}$ respectively.
In the cones $\mathring{\spc{L}}$ and $\mathring{\spc{U}}$, we use 
``$|{*}|$'' for distance to the vertex, say $o$, and 
``$\cdot$'' for cone multiplication. We also use short-cuts
$\mangle(x,y)\df\mangle\hinge{o}{x}{y}$ 
and 
\begin{align*}
 \<x,y\>\df&\,|x|\cdot|y|\cdot\cos\mangle\hinge{o}{x}{y}=\\
=&\,\tfrac12\l(|x|^2+|y|^2-\dist[2]{x}{y}{}\r).
\end{align*}
In particular,
\begin{itemize}
\item $\dist{x}{y}{\spc{L}}=\mangle(x,y)$ for any $x,y\in\spc{L}$,
\item $\dist{x}{y}{\spc{U}}=\mangle(x,y)$ for any $x,y\in\spc{U}$,
\item for any $y\in \spc{U}$, we have
$$\mangle(z,y)\le\tfrac\pi2.\eqlbl{eq:=<pi/2}$$

\end{itemize}
Set $\mathring{Q}=\Cone Q\subset \mathring{\spc{L}}$ and let $\mathring f\:\mathring{Q}\to \mathring{\spc{U}}$ be the natural cone extension of $f$; 
i.e., 
$y=f(x)$ $\Rightarrow$ $t\cdot y=\mathring f(t\cdot x)$ 
for $t\ge0$.
Clearly $\mathring f$ is short.

Applying \ref{thm:kirsz} for $\mathring f$, 
we get a short extension map $\mathring F\:\mathring{\spc{L}}\to\mathring{\spc{U}}$. 
Set $s=\mathring F(p)$.
Thus, 
$$\dist{s}{\mathring f(w)}{}
\le 
\dist{p}{w}{}
\eqlbl{eq:clm:kirszbraun-curv=1-rad-star}$$
for any $w\in \mathring Q$.
In particular, $|s|\le 1$.
Applying \ref{eq:clm:kirszbraun-curv=1-rad-star} 
for $w=t\cdot x$ and $t\to\infty$ we get

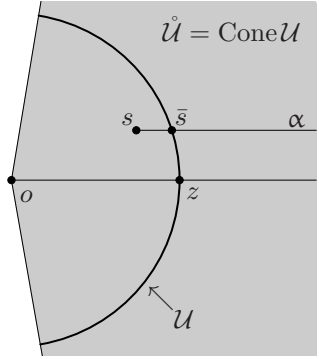
\begin{wrapfigure}[17]{l}{42mm}
\begin{lpic}[t(0mm),b(3mm),r(0mm),l(0mm)]{pics/k_0(0.3)}
\lbl[t]{100,150;$\mathring{\spc{U}}=\Cone \spc{U}$}
\lbl[br]{75,20;$\nwarrow$}
\lbl[tl]{75,20;$\spc{U}$}
\lbl[lt]{80,74;$z$}
\lbl[lb]{75,102;$\bar s$}
\lbl[b]{128,101;$\alpha$}
\lbl[rb]{57,102;$s$}
\lbl[lt]{7,73;$o$}
\end{lpic}
\end{wrapfigure}

$$\<f(x),s\>\ge \cos\mangle(p,x)\eqlbl{eq:<,>=<}$$
for any $x\in Q$.

Since $\spc{U}\in\cCat{}{0}$,
the geodesics $\geod_{[s\ t\cdot z]}$ converge as $t\to\infty$ to a ray, say $\alpha\:[0,\infty)\to \mathring{\spc{U}}$.
From \ref{eq:=<pi/2}, 
we have that the function $t\mapsto\<f(x),\alpha(t)\>$ is non-decreasing. 
Therefore, from \ref{eq:<,>=<}, for
the necessarily unique point $\bar s$ on the ray $\alpha$ such that $|\bar s|=1$ we also have 
$$\<f(x),\bar s\>\ge \cos\mangle(p,x)$$
or
$$\mangle(\bar s,f(x))
\le 
\mangle(p,f(x))$$
for any $x\in Q$.
The latter contradicts \ref{eq:cap-of-balls}.
\qeds

\section{(2\textit{n}+2)-point comparison}\label{sec:2n+2}

Here we give a generalization of the (2+2)-point comparison  to (2\textit{n}+2) points.  It follows from the generalized Kirszbraun's theorem.

First let us give a reformulation of (2+2)-point comparison. 

\begin{thm}{Reformulation of (2+2)-point comparison}
Let $\spc{X}$ be a metric space.
A quadruple $p,q,x,y\in \spc{X}$ satisfies (2+2)-point comparison if one of the following holds:
\begin{subthm}{}
One of the triples 
$(p,q,x)$ 
or 
$(p, q, y)$ 
has perimeter $>2\cdot\varpi\kappa$.
\end{subthm}

\begin{subthm}{}
If $\trig{\~p}{\~q}{\~x}
=
\modtrig\kappa(p q x)$ 
and
$\trig{\~p}{\~q}{\~y}
=
\modtrig\kappa p q y$, then
$$\dist{\~x}{\~z}{}+\dist{\~z}{\~y}{}\ge \dist{x}{y}{},$$
for any $\~z\in[\~p\~q]$.

\end{subthm}

\end{thm}

\begin{thm}{(2\textit{n}+2)-point comparison}\label{CBA-n-point}
Let $\spc{U}\in\cCat{}{\kappa}$.
Consider $x,y\in \spc{U}$ and  an array of pairs of points $(p^1,q^1)$, $(p^2,q^2),\dots,(p^n,q^n)$  in $\spc{U}$, such that there is a model configuration
$\~x$, $\~y$ and array of pairs $(\~p^1,\~q^1)$, $(\~p^2,\~q^2),\dots,(\~p^n,\~q^n)$ in $\Lob{3}\kappa$ with the following properties:
\begin{subthm}{}
$\trig{\~x}{\~p^1}{\~q^1}=\modtrig\kappa x p^1q^1$
and 
$\trig{\~y}{\~p^n}{\~q^n}=\modtrig\kappa y p^n q^n$;
\end{subthm}

\begin{subthm}{}
The simplex $\~p^i\~p^{i+1}\~q^i\~q^{i+1}$ is a model simplex%
\footnote{i.e.,
perimeter of each triple in $p^i,p^{i+1},q^i$ and $q^{i+1}$ is $<2\cdot\pi$ and
$\dist{\~p^i}{\~q^i}{}
=\dist{p^i}{q^i}{}$,
$\dist{\~p^i}{\~p^{i+1}}{}
=\dist{p^i}{p^{i+1}}{}$,
$\dist{\~q^i}{\~q^{i+1}}{}
=\dist{q^i}{q^{i+1}}{}$,
$\dist{\~p^i}{\~q^{i+1}}{}
=
\dist{p^i}{q^{i+1}}{}$ 
and $\dist{\~p^{i+1}}{\~q^{i}}{}=\dist{p^{i+1}}{q^{i}}{}$.}
 of $p^ip^{i+1}q^iq^{i+1}$
for all $i$.
\end{subthm}

Then for any choice of $n$ points $\~z^i\in [\~p^i\~q^i]$,
we have
$$\dist{\~x}{\~z^1}{}+\dist{\~z^1}{\~z^2}{}+\dots+\dist{\~z^{n-1}}{\~z^n}{}+\dist{\~z^n}{\~y}{}
\ge 
\dist{x}{y}{}.
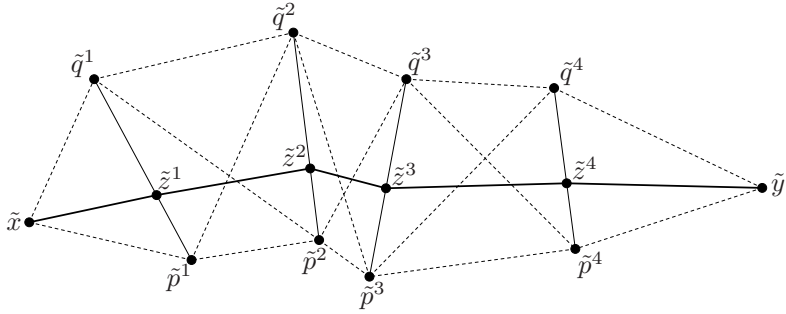
$$
\begin{center}
\begin{lpic}[t(0mm),b(0mm),r(0mm),l(0mm)]{pics/chain(0.27)}
\lbl[r]{4,33;$\~x$}
\lbl[tr]{87,12;$\~p^1$}
\lbl[t]{147,20;$\~p^2$}
\lbl[t]{175,3;$\~p^3$}
\lbl[lt]{275,18;$\~p^4$}
\lbl[br]{40,104;$\~q^1$}
\lbl[br]{138,127;$\~q^2$}
\lbl[bl]{192,105;$\~q^3$}
\lbl[bl]{266,100;$\~q^4$}
\lbl[bl]{70,49;$\~z^1$}
\lbl[br]{143,60;$\~z^2$}
\lbl[bl]{184,51;$\~z^3$}
\lbl[bl]{272,54;$\~z^4$}
\lbl[l]{369,51;$\~y$}
\end{lpic}
\end{center}
\end{thm}

To prove (2\textit{n}+2)-point comparison, we need the following lemma, which is an easy corollary from Kirszbraun's theorem (\ref{thm:kirsz+}).

\begin{thm}{Lemma}\label{cor:kir-from-hemisphere}
Let $\spc{L}\in\CBB{}{\kappa}$, 
$\spc{U}\in\cCat{}{\kappa}$,
and $Q\subset \oBall(p,\tfrac{\varpi\kappa}2)\subset \spc{L}$.
Then any short map $f\:Q\to \spc{U}$ can be extended to a short map 
$F\:\spc{L}\to \spc{U}$.
\end{thm}

\parit{Proof.} Directly from Kirszbraun's theorem (\ref{thm:kirsz} or \ref{thm:kirsz+}), we obtain the case $\kappa\le 0$. 
Thus it remains to prove the theorem only in case $\kappa>0$.
After rescaling we may assume that $\kappa=1$
and therefore $\varpi\kappa=\pi$.

It is enough to prove that there is a point $z\in \spc{U}$ such that $\dist{z}{f(x)}{}\le \tfrac\pi2$ for all $x\in Q$; once it is proved, the statement follows from Kirszbraun's theorem (\ref{thm:kirsz+}).

Further we use the same notations as in the proof of \ref{thm:kirsz+}. 

Apply Kirszbraun's theorem (\ref{thm:kirsz} or \ref{thm:kirsz+}) for $\mathring f\:\mathring Q\to\mathring{\spc{U}}$ and set $q\z={\mathring F}(p)$.
Clearly,
$$\<f(x),q\>\ge \cos\mangle(p,x)>0$$
for any $x\in Q$.
In particular, $|q|>0$. 
Thus, for $z=\tfrac{1}{|q|}\cdot q\in\spc{U}$,
we get $\dist{z}{f(x)}{\spc{U}}=\mangle(z,f(x))\le \tfrac{\pi}{2}$ for all $x\in Q$.
\qeds

\parit{Proof of (2\textit{n}+2)-point comparison.} Direct application of \ref{cor:kir-from-hemisphere} 
gives an array of short maps $f^0,f^1,\dots,f^n\:\Lob{3}\kappa\to \spc{U}$ such that
\begin{enumerate}[(i)]

\item $\~x\stackrel{f^0}{\longmapsto} x$, 
$\~p^1\stackrel{f^0}{\longmapsto} p^1$ and 
$\~q^1\stackrel{f^0}{\longmapsto}q^1$;

\item 
$\~p^i      \stackrel{f^i}{\longmapsto} p^i$, 
$\~q^{i}    \stackrel{f^i}{\longmapsto} q^i$ and 
$\~p^{i+1}  \stackrel{f^i}{\longmapsto} p^{i+1}$, 
$\~q^{i+1}  \stackrel{f^i}{\longmapsto} q^{i+1}$\\ 
for $1\le i\le n-1$;
\item 
$\~p^n\stackrel{f^n}{\longmapsto} p^n$,
$\~q^n\stackrel{f^n}{\longmapsto}q^n$ and $\~y\stackrel{f^n}{\longmapsto} y$.
\end{enumerate}
For each $i>0$, we have that $f^{i-1}|_{[\~p^i\~q^i]}=f^{i}|_{[\~p^i\~q^i]}$, since 
both $f^{i-1}$ and $f^{i}$ send $[\~p^i\~q^i]$ isometrically to a geodesic $[p^i q^i]$ in $\spc{U}$ which has to be unique
.
Thus the curves
$$f^0([\~x\~z^1]),\ f^1([\~z^1\~z^2]),\dots,\ f^{n-1}([\~z^{n-1}\~z^n]),\ f^n([\~z^n\~y])$$ 
can be joined in $\spc{U}$ into a curve connecting $x$ to $y$ with length at most 
$$\dist{\~x}{\~z^1}{}+\dist{\~z^1}{\~z^2}{}+\dots+\dist{\~z^{n-1}}{\~z^n}{}+\dist{\~z^n}{\~y}{}.\eqno\qed$$

\section{Comments and open problems}\label{sec:kirszbraun:open}

\begin{thm}{Open problem}\label{open:n-point-CBB}
Find a necessary and sufficient condition for a finite metric space to be isometrically embeddable into some $\CBB{}{\kappa}$ space.
\end{thm}

A metric on a finite set $\{a^1,a^2,\dots,a^n\}$,
can be described by the matrix with components
$$s^{ij}=\dist[2]{a^i}{a^j}{},$$
which we will call the  \emph{decrypting matrix}\index{decrypting matrix}.
The set of decrypting matrices of all metrics that admit an isometric embedding into a $\CBB{}{0}$ space 
form a convex cone, as follows from the fact that the  product of $\CBB{}{0}$ spaces is a $\CBB{}{0}$ space.
This convexity gives hope that the cone admits an explicit description.

The set of metrics on $\{a^1,a^2,\dots,a^n\}$ that can be embedded into a product of spheres with different radii admits a simple description.
Obviously, this gives a sufficient condition for \ref{open:n-point-CBB}.
This  condition is not necessary.
For instance, as follows from from a result of Vilms, \cite[2.2]{vilms},
a sufficiently dense finite subset in a generic closed positively 
curved manifold cannot be embedded into a product of spheres.

Theorem \ref{thm:pos-config} gives a necessary condition for \ref{open:n-point-CBB},
but the condition is not sufficient.
One sees this in the following example constructed by Sergei Ivanov.
A generalization of this example is given in \cite[1.1]{LPZ}.

\begin{wrapfigure}{r}{40mm}
\begin{lpic}[t(0mm),b(0mm),r(0mm),l(0mm)]{pics/ivanov-example(0.7)}
\lbl[rb]{1,16;$a$}
\lbl[lb]{55,16;$b$}
\lbl[rb]{28,31;$x$}
\lbl[rb]{28,7;$y$}
\lbl[rt]{28,1;$z$}
\lbl[lb]{44,26;$q$}
\end{lpic}
\end{wrapfigure}

\parbf{Example.}
Consider the finite set $\spc{F}\z=\{a,b,x,y,z,q\}$ with distances defined as follows:
\begin{enumerate}
\item $\dist{a}{b}{}=4$;
\item $\dist{a}{x}{}\z=\dist{a}{y}{}\z=\dist{a}{z}{}\z=\dist{b}{x}{}\z=\dist{b}{y}{}\z=\dist{b}{z}{}\z=2$;
\item $\dist{x}{y}{}=2$, $\dist{y}{z}{}=1$, $\dist{x}{z}{}=3$;
\item $\dist{x}{q}{}=\dist{q}{b}{}=1$ and thus $\dist{a}{q}{}=3$;
\item $\angk{0}{x}{q}{y}=\angk{0}{x}{q}{z}=\tfrac\pi3$; 
i.e. $\dist{q}{y}{}=\sqrt{3}$ and $\dist{q}{z}{}=\sqrt{7}$.
\end{enumerate}

On the diagram the degenerate triangles are marked by solid lines.
Note that if one removes from $\spc{F}$ the point $q$ then the remaining part can be embedded in a sphere of intrinsic diameter $4$ with poles at $a$ and $b$ and the points $x,y,z$ on the equator.
On the other hand, if one removes the  point $a$ from the space and changes the distance $\dist{z}{b}{}$ then it can be isometrically embedded into  the plane.

It is straightforward to check that this finite set satisfies the conclusion of Theorem \ref{thm:pos-config} for $\kappa=0$.
However, if such a metric appeared as an inherited metric on a subset $\{a,b,x,y,z,q\}\subset \spc{L}\in\CBB{}{0}$
then clearly  
$$\mangle\hinge x a y\z=\mangle\hinge y a z\z=\mangle\hinge y b z\z= \tfrac{\pi}{3},$$ 
contradicting $\dist{b}{z}{}=2$.

\medskip

The following problem was mentioned by Gromov in \cite[15(b)]{gromov-CAT}

\begin{thm}{Open problem}\label{open:n-point-CBA}
Describe metrics on an $n$-point set which are embeddable into $\Cat{}{\kappa}$ spaces.
\end{thm}

The set of metrics on $\{a^1,a^2,\dots,a^n\}$ which can be embedded into a product of trees and hyperbolic spaces admits a simple description using decrypting matrices defined above.
Obviously, this gives a sufficient condition for problem \ref{open:n-point-CBA}.
This  condition is not necessary.
The existence of a counterexample follows
from the same result of Vilms \cite[2.2]{vilms};
it is sufficient to take a sufficiently dense finite subset 
in a ball in a generic Hadamard space.

The (2$n$+2)-point comparison (\ref{CBA-n-point}) gives a necessary condition for \ref{open:n-point-CBA} 
which is not sufficient.
One can see this in the following example constructed by Nina Lebedeva:

Consider a square $[\~x^1\~y^1\~x^2\~y^2]$ in $\EE^3$
with  two more points  $\~z^1$, $\~z^2$ in general position on opposite sides of the plane spanned by $[\~x^1\~y^1\~x^2\~y^2]$ so that the convex hull of $\~x^1,\~x^2,\~y^1,\~y^2,\~z^1,\~z^2$ forms a nonregular octahedron with the faces formed by triangles $[\~x^i \~y^j \~z^\kay]$.
Consider the induced metric on the 6-point set $\~x^1,\~x^2,\~y^1,\~y^2,\~z^1,\~z^2$.
Note that if we increase the distance $\dist{\~z^1}{\~z^2}{}$ slightly 
then in the obtained 6-point metric $\spc{F}_6$ space all the (2+2) and (4+2)-point comparisons continue to hold.

Now assume we embed the points $x^1,x^2,y^1,y^2,z^1,z^2$ to lie in a $\cCat{}{0}$ space $\spc{U}$ in such a way that all the distances except $\dist{z^1}{z^2}{}$ are the same as between corresponding points in $\spc{F}_6$.
Note that $[\~x^1\~y^1\~x^2\~y^2]$ is a square,
therefore
we get that $\spc{U}$ contains an isometric copy of a square $\Conv(x^1,y^1,x^2,y^2)_{\spc{U}}\iso\Conv(\~x^1,\~y^1,\~x^2,\~y^2)_{\EE^3}$.
Let \[\~w\in  \Conv(\~x^1,\~y^1,\~x^2,\~y^2)_{\EE^3}\]
and $w$ be the corresponding point in $\Conv(x^1,y^1,x^2,y^2)_{\spc{U}}$.
By 
point-on-side comparison (\ref{cat-monoton}) we have $\dist{z^i}{w}{\spc{U}}\le \dist{\~z^i}{\~w}{\EE^3}$.
By construction, we may choose $\~w$ so that \[\dist{\~z^1}{\~z^2}{\EE^3}=\dist{\~z^1}{\~w}{\EE^3}+\dist{\~w}{\~z^2}{\EE^3}.\]
It follows that 
\begin{align*}
\dist{z^1}{z^2}{\spc{U}}
&\le\dist{z^1}{w}{\spc{U}}+\dist{w}{z^2}{\spc{U}}\le
\\
&\le\dist{\~z^1}{\~w}{\EE^3}+\dist{\~w}{\~z^2}{\EE^3}=
\\
&=\dist{\~z^1}{\~z^2}{\EE^3},
\end{align*}
a contradiction.

The next conjecture (if true) would give the  right generality for  Kirszbraun's theorem (\ref{thm:kirsz+}).
Roughly it states that the example \ref{example:SS_+}
is the only obstacle for extending a short map.

\begin{thm}{Conjecture}\label{conj:kirsz}
Assume $\spc{L}\in\CBB{}{1}$,
$\spc{U}\in\cCat{}{1}$,
$Q\subset \spc{L}$ is a proper subset,
and $f\: Q\to\spc{U}$ is a short map that does not admit a short extension to any bigger set $Q'\supset Q$. 
Then: 

\begin{subthm}{}
$Q$ is isometric to a sphere in a Hilbert space (of finite or cardinal dimension).
Moreover, there is a point $p\in \spc{L}$ such that $\dist{p}{q}{}=\tfrac{\pi}{2}$ for any $q\in Q$.
\end{subthm}

\begin{subthm}{}
The map $f\:Q\to\spc{U}$ is a global isometric embedding and there is no point $p'\in \spc{U}$ such that $\dist{p'}{q'}{}=\tfrac{\pi}{2}$ for all $q'\in f(Q)$.
\end{subthm}
\end{thm}

\appendix
\def\claim#1{%
\par%
\medskip%
\noindent%
\refstepcounter{thm}%
\hbox{\bf\boldmath \Alph{section}.\arabic{thm}. #1.}
\it\ 
}
\section{Barycentric simplex}\label{sec:baricentric}

The barycentric simplex was introduced by Kleiner in \cite{kleiner};
it is a construction that works in a general metric space.
Roughly, it gives a $\kay$-\nospace dimensional submanifold for a given ``nondegenerate'' array of $\kay+1$ strongly convex functions.

Let us denote by $\Delta^\kay\subset \RR^{\kay+1}$\index{$\Delta^m$} 
the \emph{standard $\kay$-simplex}\index{standard simplex}; 
i.e. $\bm{x}=(x_0,x_1,\dots,x_n)\in\Delta^\kay$ if $\sum_{i=0}^\kay x_i=1$ and $x_i\ge0$ for all $i$.

Let $\spc{X}$ be a metric space 
and $\bm{f}=(f^0,f^1,\dots,f^\kay)\:\spc{X}\to \RR^{\kay+1}$ be a function array.
Consider the map $\spx{\bm{f}}\:\Delta^\kay\to \spc{X}$,\index{$\spx{\bm{f}}$} defined by 
$$\spx{\bm{f}}(\bm{x})=\argmin\sum_{i=0}^\kay x_i\cdot f^i,$$
where $\argmin f$\index{$\argmin$} denotes a point of minimum of $f$.
The map $\spx{\bm{f}}$ will be called a \emph{barycentric simplex}\index{barycentric simplex of function array} of $\bm{f}$.
In general, a barycentric simplex of a function array might be undefined and need not be unique. 

The name comes from the fact that if $\spc{X}$ is a Euclidean space 
and $f^i(x)\z=\tfrac{1}{2}\cdot\dist[2]{p^i}{x}{}$ for some array of points $\bm{p}=(p^0,p^1,\dots,p^\kay)$, 
then $\spx{\bm{f}}(\bm{x})$ is the barycenter of points $p^i$ with weights $x_i$.
 
A barycentric simplex $\spx{\bm{f}}$ 
for the function array $f^i(x)=\tfrac{1}{2}\cdot\dist[2]{p^i}{x}{}$
will also be called a \emph{barycentric simplex with vertices at} $\{p^i\}$\index{barycentric simplex with vertices at point array}.

It is clear from the  definition that if 
$\bm{\hat f}$ is a subarray of $\bm{f}$,
then $\spx{\bm{\hat f}}$ coincides with the restriction of $\spx{\bm{f}}$ to the corresponding face of $\Delta^\kay$.

The following theorem shows that the barycentric simplex is defined 
for an array of strongly convex functions on a complete geodesic space. 
In order to formulate the theorem, we need to introduce a partial order $\succcurlyeq$ on $\RR^{\kay+1}$.

\begin{thm}{Definition}\label{def:supset+succcurlyeq}
For two real arrays $\bm{v}$, $\bm{w}\in \RR^{\kay+1}$,
$\bm{v}=(v^0,v^1,\dots,v^\kay)$ 
and 
$\bm{w}=(w^0,w^1,\dots,w^\kay)$, 
we write
$\bm{v}\succcurlyeq\bm{w}$ if $v^i\ge w^i$ for each $i$.

Given a subset $Q\subset \RR^{\kay+1}$, define its \emph{superset}\index{superset}
\index{$\SupSet$}
$$\SupSet Q =\{\bm{v}\in\RR^\kay\mid\exists\, \bm{w}\in Q\ \t{such that}\ \bm{v}\succcurlyeq\bm{w}\}.$$

\end{thm}

\begin{thm}{Theorem on barycentric simplex}\label{thm:bary}
Assume $\spc{X}$ is a complete geodesic space and 
$\bm{f}\z=(f^0,f^1,\dots,f^\kay)\:\spc{X}\to\RR^\kay$ is an array of strongly convex and locally Lipschitz functions.

Then the barycentric simplex $\spx{\bm{f}}\:\Delta^\kay\to \spc{X}$
is uniquely defined and moreover:

\begin{subthm}{bary-Lip} $\spx{\bm{f}}$ is Lipschitz. 
\end{subthm}

\begin{subthm}{bary-iff} The set $\SupSet{\bm{f}(\spc{X})}\subset\RR^{\kay+1}$ is convex,
and
$p\in \spx{\bm{f}}(\Delta^\kay)$ if and only if
$\bm{f}(p)\in\Fr\l[\SupSet{\bm{f}(\spc{X})}\r]$.
In particular, $\bm{f}\circ\spx{\bm{f}}(\Delta^\kay)$ lies on a convex hypersurface in $\RR^{\kay+1}$.
\end{subthm}

\begin{subthm}{bary-embed} The restriction $\bm{f}|_{\spx{\bm{f}}(\Delta^\kay)}$  has  $C^{\frac{1}{2}}$-inverse.
\end{subthm}

\begin{subthm}{bary-R^n} 
The set $\mathfrak{S}=\spx{\bm{f}}(\Delta^\kay)\backslash\spx{\bm{f}}(\partial\Delta^\kay)$
is $C^{\frac{1}{2}}$-homeomorphic to an open domain in $\RR^\kay$.
\end{subthm}
\end{thm}

The set $\mathfrak{S}$ described above will be called \emph{Kleiner's spine}\index{Kleiner's spine} of $\bm{f}$.
If $\mathfrak{S}$ is nonempty, we say the barycentric simplex $\spx{\bm{f}}$ is \emph{nondegenerate}\index{nondegenerate}.

We precede the proof of the theorem with the following lemma.

\begin{thm}{Lemma}\label{lem:argmin(convex)}
Assume $\spc{X}$ is a complete geodesic metric space and let  $f\:\spc{X}\to\RR$ be a locally Lipschitz, strongly convex function.  Then the minimum point 
$p=\argmin f$ 
is uniquely defined.
\end{thm}

\parit{Proof.}
Assume that $x$ and $y$ are distinct minimum points of $f$. 
Then for the midpoint $z$ of a geodesic $[x y]$ we have
$$f(z)<f(x)=f(y),$$ 
a contradiction. 
It only remains to show existence.

Fix a point $p\in  \spc{X}$; 
let $\Lip\in\RR$ be a Lipschitz constant of $f$ in a neighborhood of $p$.
Without loss of generality, we can assume that $f$ is $1$-convex.
Consider function $\phi(t)=f\circ\geod_{[px]}(t)$.
Clearly $\phi$ is $1$-convex and $\phi^+(0)\ge -\Lip$.
Setting $\ell=\dist{p}{x}{}$, we get 
\begin{align*}
f(x)
&=
\phi(\ell)
\ge
\\
&\ge
f(p)-\Lip\cdot\ell+\tfrac{1}{2}\cdot\ell^2
\ge
\\
&\ge f(p)-\tfrac{1}{2}\cdot{\Lip^2}.
\end{align*}

In particular,
$$s
\df
\inf\set{f(x)}{x\in \spc{X}}
\ge
f(p)-\tfrac{1}{2}\cdot{\Lip^2}.$$
If $z$ is a midpoint of $[x y]$ then  
$$s\le f(z)
\le
\tfrac{1}{2}\cdot f(x)+\tfrac{1}{2}\cdot f(y)-\tfrac{1}{8}\cdot\dist[2]{x}{y}{}.
\eqlbl{mid-point}$$
Choose a sequence of points $p_n\in \spc{X}$  such that $f(p_n)\to s$.
Applying \ref{mid-point}, for $x\z=p_n$, $y\z=p_m$, we get that $(p_n)$ is a Cauchy sequence. 
Clearly, $p_n\to \argmin f$.
\qeds

\parit{Proof of theorem \ref{thm:bary}.}
Without loss of generality, we can assume that each $f^i$ is $1$-convex.
Thus, for any $\bm{x}\in\Delta^\kay$, 
the convex combination $\sum x_i\cdot f^i\:\spc{X}\to\RR$ is also $1$-convex.
Therefore, according to Lemma~\ref{lem:argmin(convex)}, $\spx{\bm{f}}(\bm{x})$ is defined.

\parit{(\ref{SHORT.bary-Lip}).} 
Since $\Delta^\kay$ is compact, it is sufficient to show that $\spx{\bm{f}}$ is locally Lipschitz.

For $\bm{x},\bm{y}\in\Delta^\kay$,
set 
\begin{align*}
f_{\bm{x}}
&=\sum x_i\cdot f^i,
&
f_{\bm{y}}
&=\sum y_i\cdot f^i,
\\
p
&=\spx{\bm{f}}(\bm{x}),
&
q
&=\spx{\bm{f}}(\bm{y}).
\end{align*}
Let $\ell=\dist[2]{p}{q}{}$.
Clearly 
$\phi(t)=f_{\bm{x}}\circ\geod_{[p q]}(t)$ takes its minimum at $0$ and
$\psi(t)=f_{\bm{y}}\circ\geod_{[p q]}(t)$ takes its minimum at $\ell$.
Thus $\phi^+(0)$, $\psi^-(\ell)\ge 0$%
\footnote{Here $\phi^\pm$ denotes ``signed one sided derivative''; i.e. 
$$\phi^\pm(t_0)=\lim_{t\to t_0\pm}\frac{\phi(t)-\phi(t_0)}{|t-t_0|}$$}%
.
From $1$-convexity of $f_{\bm{y}}$, we have
$\psi^+(0)+\psi^-(\ell)+\ell\le0$.

Let $\Lip$ be a Lipschitz constant for all $f^i$ in a neighborhood $\Omega\ni p$.
Then $\psi^+(0)\le \phi^+(0)+\Lip\cdot\|\bm{x}-\bm{y}\|_{{}_1}$, 
where $\|\bm{x}-\bm{y}\|_{{}_1}=\sum_{i=0}^\kay|x_i-y_i|$.
That is, given $\bm{x}\in\Delta^\kay$, there is a constant $\Lip$ such that
$$\dist{\spx{\bm{f}}(\bm{x})}{\spx{\bm{f}}(\bm{y})}{}
=
\ell\le \Lip\cdot\|\bm{x}-\bm{y}\|_{{}_1}$$
for any $\bm{y}\in\Delta^\kay$.
In particular, there is $\eps>0$ such that if $\|\bm{x}-\bm{y}\|_{{}_1},$ $\|\bm{x}-\bm{z}\|_{{}_1} <\eps$, then $\spx{\bm{f}}(\bm{y})$, $\spx{\bm{f}}(\bm{z})\in\Omega$. 
Thus, the same argument as above implies 
$$\dist{\spx{\bm{f}}(\bm{y})}{\spx{\bm{f}}(\bm{z})}{}
=
\ell\le \Lip\cdot\|\bm{y}-\bm{z}\|_{{}_1}$$
for any $\bm{y}$ and $\bm{z}$ sufficiently close to $\bm{x}$;
i.e. $\spx{\bm{f}}$ is locally Lipschitz.

\parit{(\ref{SHORT.bary-iff}).} The ``only if'' part is trivial, let us prove  the ``if''-part.

Note that convexity of $f^i$ implies that
for any two points $p,q\in \spc{X}$ and $t\in[0,1]$ we have
$$(1-t)\cdot\bm{f}(p)+t\cdot \bm{f}(q)
\succcurlyeq
\bm{f}\circ\geodpath_{[p q]}(t),
\eqlbl{n-convex}$$
where $\geodpath_{[p q]}$ is a geodesic path from $p$ to $q$; 
i.e. $\geodpath_{[p q]}(t)=\geod_{[p q]}(\tfrac{t}{\dist{p}{q}{}})$. 

From \ref{n-convex}, we have that $\SupSet[\bm{f}(\spc{X})]$ is a convex subset of $\RR^{\kay+1}$.
If 
$$\max_{i}\{f^i(q)\z-f^i(p)\}\ge0$$ 
for any $q\in \spc{X}$, then $\bm{f}(p)$ lies in the boundary of $\SupSet[\bm{f}(\spc{X})]$.
Take a supporting vector $\bm{x}\in\RR^{\kay+1}$ to $\SupSet[\bm{f}(\spc{X})]$ at $\bm{f}(p)$.
Thus $\bm{x}\not=\bm{0}$ and $\sum_i x_i\cdot[w^i-f^i(p)]\ge0$ for any 
$\bm{w}\in \SupSet[\bm{f}(\spc{X})]$. In particular, $\sum_i x_i\cdot v_i \ge 0$ for any $v=(v_1,\ldots, v_k)$ with all $v_i\ge 0$.  Hence $x_i\ge 0$ for all $i$ and 
$\bm{x}'=\frac{\bm{x}}{\|\bm{x}\|}_{{}_1}\in\Delta^\kay$. 
Thus $p=\spx{\bm{f}}(\bm{x}')$. 

\parit{(\ref{SHORT.bary-embed}).}
The restriction $\bm{f}|_{\spx{\bm{f}}(\Delta^\kay)}$ is Lipschitz.
Thus we only have to show that it has a  $C^{\frac{1}{2}}$-inverse.
Given $\bm{v}\in\RR^{\kay+1}$, consider the function 
$h_{\bm{v}}\: \spc{X}\to \RR$ given by
$$h_{\bm{v}}(p)=\max_{i}\{f^i(p)-v^i\}.$$
Define a map $\map \:\RR^{\kay+1}\to \spc{X}$ by
$\map (\bm{v})=\argmin h_{\bm{v}}$.

Clearly $h_{\bm{v}}$ is $1$-convex.
Thus, according to \ref{lem:argmin(convex)}, $\map (\bm{v})$ is uniquely defined for any $\bm{v}\in\RR^{\kay+1}$.
From (\ref{SHORT.bary-iff}), for any $p\in \spx{\bm{f}}(\Delta^\kay)$ we have
$\map \circ\bm{f}(p)=p$.

It remains to show that $\map $ is $C^{\frac{1}{2}}$-continuous.
Clearly, 
$$|h_{\bm{v}}-h_{\bm{w}}|
\le\|\bm{v}-\bm{w}\|_\subinfty
\df
\max_{i}\{|v^i-w^i|\},$$
for any $\bm{v},\bm{w}\in\RR^{\kay+1}$.
Set $p=\map (\bm{v})$ and $q=\map (\bm{w})$.
Since $h_{\bm{v}}$ and $h_{\bm{w}}$ are 1-convex,
\begin{align*}
h_{\bm{v}}(q)
&\ge 
h_{\bm{v}}(p)+\tfrac{1}{2}\cdot\dist[2]{p}{q}{},
&
h_{\bm{w}}(p)
&\ge 
h_{\bm{w}}(q)+\tfrac{1}{2}\cdot\dist[2]{p}{q}{}.
\end{align*}
Therefore,
$$\dist[2]{p}{q}{}\le 2\cdot\|\bm{v}-\bm{w}\|_\subinfty.$$
Hence the result.

\parit{(\ref{SHORT.bary-R^n}).} 
Let $S=\Fr\SupSet(\bm{f}(\spc{X}))$.
Note that orthogonal projection to the hyperplane $\WW^\kay$ in $\RR^{\kay+1}$ defined by equation $x_0+x_1+\dots+x_n=0$ gives a bi-Lipschits homeomorphism $S\to \WW^\kay$.

Clearly, $\bm{f}({\spx{\bm{f}}(\Delta^\kay)}\backslash\spx{\bm{f}}(\partial\Delta^\kay))$ 
is an open subset of $S$.
Hence the result.
\qeds
\section{Helly's theorem}\label{sec:helly}

\begin{thm}{Helly's theorem}\label{thm:helly}
Let $\spc{U}\in\cCat{}0$ 
and $\{K_\alpha\}_{\alpha\in \IndexSet}$ be an arbitrary collection of closed bounded convex subsets of $\spc{U}$.

If 
$$\bigcap_{\alpha\in \IndexSet}K_\alpha=\emptyset,$$
then there is an index array $\alpha_1,\alpha_2,\dots,\alpha_n\in \IndexSet$ such that
$$\bigcap_{i=1}^nK_{\alpha_i}=\emptyset.$$

\end{thm}

\parbf{Remarks.}
\begin{enumerate}[(i)]
\item In general, none of $K_\alpha$ might be compact. 
Thus the the statement is not completely trivial.
\item If $\spc{U}$ is a Hilbert space (not necessarily separable), 
then the above result is equivalent to the statement that a convex bounded set 
which is closed in the ordinary topology forms a compact set in the weak topology.

In fact, one can define the \emph{weak topology} on an arbitrary metric space, by taking exteriors of closed balls as its prebase.
Then the result above implies for $\spc{U}\in\cCat{}0$, any closed bounded convex set in $\spc{U}$ is compact in the  weak topology 
(this is very similar to the definition given by Monod in \cite{monod}).
\end{enumerate}

We present the proof of Lang and Shroeder from \cite{lang-schroeder}.

\begin{thm}{Lemma}\label{lem:closest point}
Let $\spc{U}\in\cCat{}0$.
Given a closed convex set $K\subset \spc{U}$ and a point $p\in \spc{U}\backslash K$, 
there is unique point $p^*\in K$ such that $\dist{p^*}{p}{}=\dist{K}{p}{}$. 
\end{thm}

\parit{Proof.}
Let us first prove uniqueness. 
Assume there are two points $y',y''\in K$ 
so that $\dist{y'}{p}{}=\dist{y''}{p}{}=\dist{K}{p}{}$.
Take $z$ to be the midpoint of $[y'y'']$. 
Since $K$ is convex, $z\in K$.
From comparison, we have that $\dist{z}{p}{}<\dist{y'}{p}{}=\dist{K}{p}{}$, a contradiction

The proof of existence is analogous.
Take a sequence  of points $y_n\in K$ 
such that $\dist{y_n}{p}{}\to \dist{K}{p}{}$.
It is enough to show that $(y_n)$ is a Cauchy sequence; 
thus one could take $p^*=\lim_n y_n$.

Assume $(y_n)$ is not Cauchy, then for some fixed $\eps>0$, 
we can choose two subsequences $(y_n')$ and $(y_n'')$ of $(y_n)$ 
such that 
$\dist{y'_n}{y''_n}{}\ge\eps$ for each $n$.
Set $z_n$ to be the midpoint of $[y'_ny''_n]$; from convexity we have $z_n\in K$.
From point-on-side comparison (see page \pageref{POS-CAT}), there is $\delta>0$ 
such that $\dist{p}{z_n}{}\le \max\{\dist{p}{y'_n}{},\dist{p}{y''_n}{}\}-\delta$. 
Thus 
$$\limsup_{n\to\infty}\dist{p}{z_n}{}<\dist{K}{x}{},$$ 
a contradiction\qeds

\parit{Proof of \ref{thm:helly}.} 
Assume the contrary. Then for any finite set $F\subset \IndexSet$,
$$K_{F}\df \bigcap_{\alpha\in F}K_{\alpha}\not=\emptyset.$$
We construct a point $z$ such that $z\in K_\alpha$ for each $\alpha\in \IndexSet$.
Thus we arrive at a contradiction since
$$\bigcap_{\alpha\in \IndexSet}K_\alpha=\emptyset.$$

Choose a point $p\in \spc{U}$ and set $r=\sup\dist{K_{F}}{p}{}$ where $F$ runs over all finite subsets of $\IndexSet$.
Let $p^*_F$  be the closest point on $K_{F}$ from $p$; 
according to Lemma \ref{lem:closest point}, $p^*_F$ 
exits and is unique.

Take a nested sequence of finite subsets 
$F_1\subset F_2\subset \dots$ of $\IndexSet$, such that $\dist{K_{F_n}}{p}{}\to r$.

Let us show that $p^*_{F_n}$ is a Cauchy sequence. 
Indeed, if not then for some fixed $\eps>0$, 
we can choose two subsequences $(y'_n)$ and $(y''_n)$ of $(p^*_{F_n})$ 
such that $\dist{y'_n}{y''_n}{}\ge\eps$.
Set $z_n$ to be midpoint of $[y'_ny''_n]$. 
From point-on-side comparison (see page \pageref{POS-CAT}), 
there is $\delta>0$ such that 
$\dist{p}{z_n}{}\le \max\{\dist{p}{y'_n}{},\dist{p}{y''_n}{}\}-\delta$.
Thus 
$$\limsup_{n\to\infty}\dist{p}{z_n}{}<r.$$
On the other hand, from convexity, each $F_n$ 
contains all $z_\kay$ with sufficiently large $\kay$, a contradiction.

Thus, $p^*_{F_n}$ converges and we can set $z=\lim_n p^*_{F_n}$.
Clearly $\dist{p}{z}{}=r$.

Repeat the above arguments for  the sequence $F_n'=F_n\cup \{\alpha\}$.
As a result, we get another point $z'$ such that $\dist{p}{z}{}=\dist{p}{z'}{}=r$ and 
$z,z'\in K_{F_n}$ for all $n$.
Thus, if $z\not=z'$ the midpoint $\hat z$ of $[zz']$ would belong to all 
$K_{F_n}$ and from comparison we would have $\dist{p}{\hat z}{}<r$, a contradiction.

Thus, $z'=z$; in particular 
$z\in K_\alpha$ for each $\alpha\in\IndexSet$.
\qeds

\end{document}